\documentclass[12pt]{amsart}
\usepackage{mathrsfs}

\usepackage{amssymb,amsmath,graphicx,color,textcomp, amsthm,bbm,bbold, enumerate,booktabs}
\usepackage{chngcntr}
\usepackage{apptools}
\AtAppendix{\counterwithin{theorem}{section}}
\definecolor{Red}{cmyk}{0,1,1,0}

\definecolor{verde}{cmyk}{1,0,1,0}

\definecolor{loka}{cmyk}{.5,0,1,.5}
\definecolor{azul}{cmyk}{1,1,0,0}

\numberwithin{equation}{section}

\newcommand{\be}{\begin{equation}}
\newcommand{\ee}{\end{equation}}

\newtheorem{definition}{Definition}

\begin{document}
\title{A note on the Laplace transform and the variable-order differential operators}
\author{J. Vanterler da C. Sousa$^1$}
\address{$^1$ Department of Applied Mathematics, Institute of Mathematics,
 Statistics and Scientific Computation, University of Campinas --
UNICAMP, rua S\'ergio Buarque de Holanda 651,
13083--859, Campinas SP, Brazil\newline
e-mail: {\itshape \texttt{ra160908@ime.unicamp.br, capelas@ime.unicamp.br }}}

\author{E. Capelas de Oliveira$^1$}

\begin{abstract} 
In this short note, using the variable-order differential operator introduced by means of the inverse Laplace transform \cite{coimbra}, we questioned the result obtained by Yang and Tenreiro Machado \cite{yang}.    

\vskip.5cm
\noindent
\emph{Keywords}: Fractional derivative of variable-order; Laplace transform.
\newline 
MSC 2010 subject classifications. 26A33; 44A10.
\end{abstract}
\maketitle

%%%%%%%%%%%%%%%%%%%%%%%%%%%%%%%%%%%%%%%%%%%%%%%%%%%%%%%%%%%%%%%%%%%%%%%%%%%%%%%%%%%%%%%%%%%%%%%%%%%%%%%%%%%%%%%%%%%%%%%%%%%%%%%%%%%%%%%%%%%%%%%%%%%%%%%%%%%%%

\section{Introduction} 

Over the time, fractional calculus has provided researchers with several important discoveries, in particular, 
applications in various areas of knowledge \cite{AHMJ,SAM,ATA}. Some authors have shown themselves to be 
interested in proposing new operators of fractional differentiation and integration, until then with non-variable order. 
Recently, Almeida \cite{RCA} introduced the $\psi$-Caputo fractional derivative of non-variable order of a function 
$f$ with respect to another function $\varphi$ and, Yang and Tenreiro Machado \cite{yang}, presented a new operator 
of fractional integration of variable order and consequently made use of it, to introduce an operator, the $\psi$-Caputo operator of 
variable order with respect to another function.

Fractional differential operators in the Caputo or Riemann-Liouville sense, among others, to model certain natural 
phenomena, has been one of the most important tool to work with fractional calculus \cite{AHMJ,ATA}. The Laplace transform 
is widely used to solve such problems. In this note, we present a discussion about the Laplace transform of the Caputo fractional 
derivative of variable order with respect to another function, recently proposed by Yang and Tenreiro Machado \cite{yang}. 
By means of the variable order differential operator introduced as an inverse Laplace transform \cite{coimbra}, we present 
some results, specifically, we evaluate the Laplace transform of the fractional derivative and make some comparisons.

%%%%%%%%%%%%%%%%%%%%%%%%%%%%%%%%%%%%%%%%%%%%%%%%%%%%%%%%%%%%%%%%%%%%%%%%%%%%%%%%%%%%%%%%%%%%%%%%%%%%%%%%%%%%%%%%%%%%%%%%%%%%%%%%%%%%%%%%%%%%%%%%%%%%%%%%%%%%%%%%%%%%%%%%%%%%%%%%%%%%%%%%%%%%%%%%%%%%%%%%%%%%
\section{Preliminaries}
In this section we recover some results which are useful in the sequel of the note. 
First, we review the Laplace transform and the corresponding inverse Laplace transform; second we introduce the 
$\varphi $-Caputo fractional derivative of the function $\psi \left( \sigma ,t\right) $ of two-variable order 
$\xi \left( \sigma ,t\right) $ $\left( 0<\xi \left( \sigma ,t\right) <1\right) $ with respect to another function $\varphi $. 
Then, we take the Laplace transform of the $\varphi $-Caputo fractional derivative.

\begin{definition}
Let $\psi(t)$ be a function of exponential order. The Laplace transform operator, denote by $\mathscr{L}(\cdot)$, acting on a 
function $\psi(t)$, defined by means of an improper integral with a non singular kernel, is
\begin{equation}\label{1}
\mathscr{L}[\psi(t)] := \Psi(s) = \int_0^{\infty} e^{-st} \psi(t) \, {\rm{d}}t , \qquad {\mbox{Re}}(s)>0
\end{equation}
with $s \in \mathbb{C}$ is the parameter associated with the Laplace transform.
\end{definition}

The corresponding inverse Laplace transform is obtained by means of the following expression
\begin{equation}\label{2}
\mathscr{L}^{-1}[\Psi(s)]:=\psi(t)= \frac{1}{2\pi i} \int_{c-i\infty}^{c+i\infty} e^{st} \Psi(s) \, {\mbox{d}}s
\end{equation}
with $t > 0$ and the contour, in the complex plane, contains all singularities on the left of the straightright line 
${\mbox{Re}}(s)=c > 0$.

Suppose that $T$ is an interval $0\leq a<b<\infty $, $\varphi \in C^{1}(T)$ with $\varphi^{\left( 1\right) }\left( t\right) \neq 0$ 
and $\psi \in C^{1}\left( T\right) $, for $\forall t\in T.$ The left and right fractional integral of the function 
$\psi \left( \sigma ,t\right) $ of two-variable order $\xi \left( \sigma ,t\right) $ $\left( 0<\xi \left( \sigma ,t\right) <1\right),$ 
with respect to another function $\varphi $ are given by \cite{yang}
\begin{equation*}
I_{a+}^{\xi \left( \sigma ,t\right) ;\varphi }\psi \left( \sigma ,t\right) :=%
\frac{1}{\Gamma \left( \xi \left( \sigma ,t\right) \right) }%
\int_{a}^{t}\varphi ^{\left( 1\right) }\left( \sigma ,s\right) \left(
\varphi \left( \sigma ,t\right) -\varphi \left( \sigma ,s\right) \right)
^{\xi \left( \sigma ,t\right) -1}\psi \left( \sigma ,s\right) {\rm{d}}s
\end{equation*}%
and 
\begin{equation*}
I_{b-}^{\xi \left( \sigma ,t\right) ;\varphi }\psi \left( \sigma ,t\right) :=%
\frac{1}{\Gamma \left( \xi \left( \sigma ,t\right) \right) }%
\int_{t}^{b}\varphi ^{\left( 1\right) }\left( \sigma ,s\right) \left(
\varphi \left( \sigma ,s\right) -\varphi \left( \sigma ,t\right) \right)
^{\xi \left( \sigma ,t\right) -1}\psi \left( \sigma ,s\right) {\rm{d}}s,
\end{equation*}
respectively, where $\Gamma \left( \cdot \right) $ is the gamma function.

The left and right $\varphi$-Caputo fractional derivatives of the function $\psi \left( \sigma ,t\right)$ with two-variable order 
$\xi \left( \sigma ,t\right) $ $\left( 0<\xi \left( \sigma ,t\right) <1\right)$, with respect to another function $\varphi$, 
are defined by \cite{yang}
\begin{equation}\label{jose}
^{C}\mathbb{D}_{a+}^{\xi \left( \sigma ,t\right) ;\varphi }\psi \left( \sigma
,t\right) :=\frac{1}{\Gamma \left( 1-\xi \left( \sigma ,t\right) \right) } \int_{a}^{t}\left( \varphi \left( \sigma ,t\right) -
\varphi \left( \sigma ,s\right) \right) ^{-\xi \left( \sigma ,t\right) }\psi ^{\left( 1\right) }\left( \sigma ,s\right) {\rm{d}}s
\end{equation}
and
\begin{equation}\label{jose1}
^{C}\mathbb{D}_{b-}^{\xi \left( \sigma ,t\right) ;\varphi }\psi \left( \sigma
,t\right) :=\frac{-1}{\Gamma \left( 1-\xi \left( \sigma ,t\right) \right) } \int_{t}^{b}\left( \varphi \left( \sigma ,s\right) -
\varphi \left( \sigma,t\right) \right) ^{-\xi \left( \sigma ,t\right) }\psi ^{\left( 1\right) }\left( \sigma ,s\right) {\rm{d}}s,
\end{equation}
respectively, and $\psi^{\left( 1\right) }\left( \cdot,\cdot\right) $ denotes the first order derivative.

Introducing $\xi(\sigma,t)=\xi$ and $\varphi(\sigma,t)=t$ in Eq.(\ref{jose}) and Eq.(\ref{jose1}), we have the left and right 
Caputo fractional derivatives
\begin{equation}\label{jose2}
^{C}\mathbb{D}_{a+}^{\xi}\psi \left( t\right) :=\frac{1}{\Gamma \left(
1-\xi \right) }\int_{a}^{t}\left( t-s\right) ^{-\xi }\psi ^{\left( 1\right)
}\left( s\right) {\rm{d}}s
\end{equation}
and
\begin{equation}\label{jose3}
^{C}\mathbb{D}_{b-}^{\xi}\psi \left( t\right) :=\frac{-1}{\Gamma \left(
1-\xi \right) }\int_{t}^{b}\left( s-t\right) ^{-\xi }\psi ^{\left( 1\right)
}\left( s\right) {\rm{d}}s,
\end{equation}
respectively, with $0<\xi <1$ and $\psi^{\left( 1\right) }\left( \cdot\right)$ denotes the first order derivative.

First, we consider $a=0$ in Eq.(\ref{jose2}) and taking the respective Laplace transform, we get
\begin{equation}\label{4}
\mathscr{L}[{}^C {\mathbb{D}}_{0^{+}}^{\xi} \psi(t)] = s^{\xi} \Psi(s) - \sum_{k=0}^{n-1} s^{n-1-k} \psi^{(k)}(0).
\end{equation}

Considering $n=1$ in Eq.(\ref{4}) we have $\mathscr{L}[{}^C{\mathbb{D}}_{0^{+}}^{\xi} (f(t))] = s^{\xi}\Psi(s) - \psi(0)$.

For continuity of work it is appropriate to consider the following function: $\psi(t) = \dfrac{t^{-1-\alpha}}{\Gamma(-\alpha)}$, 
where $t$ is the independent variable and $\alpha$, with $0 < \alpha < 1$, is a real constant. Taking its Laplace transform, we have 
\begin{equation}  \label{LT}
\mathscr{L}[\psi(t)] = \mathscr{L} \left[\frac{t^{-1-\alpha}}{\Gamma(-\alpha)}\right] = 
\frac{1}{\Gamma(-\alpha)} \mathscr{L}[t^{-1-\alpha}].
\end{equation}

\section{Comments on the paper by Yang \& Tenreiro Machado}
In a recent paper \cite{coimbra} the author introduces a differential operator of variable order by means of the inverse Laplace 
transform of the function $s^{q(t^{\prime })}$, defined by \cite{yang}
\begin{equation}  \label{Coimbra}
\mathscr{L}^{-1} [s^{q(t^{\prime })}] = \frac{t^{-1-q(t^{\prime })}}{%
\Gamma(-q(t^{\prime }))}.
\end{equation}

It is important to note that: the variable associated with the inverse Laplace transform is $s$, with $s \neq t^{\prime }$ 
and $q(t^{\prime })$ is the order of the derivative of variable-order. Then, taking the Laplace transform on both sides 
of Eq.(\ref{Coimbra}) we get 
\begin{equation}  \label{LT1}
s^{q(t^{\prime })} = \mathscr{L} \left[ \frac{t^{-1-q(t^{\prime })}}{%
\Gamma(-q(t^{\prime }))} \right] = \frac{1}{\Gamma(-q(t^{\prime }))} %
\mathscr{L}[t^{-1-q(t')}]
\end{equation}
because $t$ is considered the variable of integration in the Laplace transform.

Choosing $q(t^{\prime })=q$, a constant in Eq.(\ref{LT1}), we have
\begin{equation*}
s^{q}=\frac{1}{\Gamma (-\alpha )}\mathscr{L}[t^{-1-q}],
\end{equation*}
which is the same as in Eq.(\ref{LT}) for $q=\alpha $.

On the other hand, as defined in \cite{coimbra} we conclude  
\begin{equation*}
\mathscr{L}\left[ \frac{t^{-1-q(t^{\prime })}}{\Gamma (-q(t^{\prime }))}\right] =
\frac{1}{\Gamma (-q(t^{\prime }))}\mathscr{L} [t^{-1-q(t')}]\neq \mathscr{L}\left[ \frac{t^{-1-q(t)}}{\Gamma (-q(t))}\right] 
\end{equation*}
because $t\neq t^{\prime }$.

Also, in this sense, in a recent paper \cite{yang}, the authors write its Eq.(15) in the following form
\begin{equation*}
\mathscr{L}\left[ \frac{t^{-1-\xi (t)}}{\Gamma (-\xi (t))}\right] =s^{\xi (t)}.
\end{equation*}

But, as we have just shown it should be written as follows
\begin{equation*}
\mathscr{L}\left[ \frac{t^{-1-\xi (t^{\prime })}}{\Gamma (-\xi (t^{\prime }))
}\right] =s^{\xi (t^{\prime })}
\end{equation*}
or in the following form $\displaystyle \frac{1}{\Gamma (-\xi (t^{\prime }))}\mathscr{L}[t^{-1-\xi(t')}]=s^{\xi (t^{\prime })}$.

The authors go beyond, in the same paper \cite{yang}, they consider $\displaystyle \dfrac{t^{-1-\xi(\sigma ,
t^{\prime })}}{\Gamma (-\xi (\sigma ,t^{\prime }))}$ with $\sigma $ a real constant, and conclude that (its Eq.(16)) 
\begin{equation}\label{jose5}
\frac{1}{\Gamma (-\xi (\sigma ,t^{\prime }))}\mathscr{L}[t^{-1-\xi(%
\sigma,t')}]=s^{\xi (\sigma ,t^{\prime })}
\end{equation}
with $t\neq \sigma \neq t^{\prime }$.

So, if we choose $\sigma =0$ in the last equation, we get its Eq.(15), with the identification 
$\xi(0,t') = \xi(t')$ and, also, if $\sigma = t'=0$ we recover Eq.(\ref{LT}) with the identification 
$\xi(0,0)=\alpha$.

Now, in the same paper \cite{yang}, the authors consider a function $\varphi (t)$, defined by 
\begin{equation*}
\frac{(\varphi (t))^{-1-\xi (\sigma ,t^{\prime })}}{\Gamma \left( -\xi
(\sigma ,t^{\prime })\right) }
\end{equation*}
whose Laplace transform was given as 
\begin{equation*}
\mathscr{L}\left[ \frac{(\varphi (t))^{-1-\xi (\sigma ,t^{\prime })}}{\Gamma
\left( -\xi (\sigma ,t^{\prime })\right) }\right] =\left[ \varphi
^{(1)}\left( t\right) \right] s^{\xi (\sigma ,t^{\prime })}
\end{equation*}%
which can be written in the following form
\begin{equation}\label{jose6}
\frac{1}{\Gamma (-\xi (\sigma ,t^{\prime }))}\mathscr{L}\left[ (\varphi
(t))^{-1-\xi (\sigma ,t^{\prime })}\right] =\left[ \varphi ^{(1)}\left(
t\right) \right] s^{\xi (\sigma ,t^{\prime })}.
\end{equation}

For $\varphi (t)=t$, we have $\varphi ^{\left( 1\right) }(t)=1$. So that Eq.(\ref{jose6}) becomes Eq.(\ref{jose5}), that is
\begin{equation}
\left[ \varphi ^{\left( 1\right) }(t)\right] s^{\xi \left( \sigma ,t^{\prime }\right) }=\left[ 1\right] 
s^{\xi \left( \sigma ,t^{\prime }\right) }=s^{\xi \left( \sigma ,t^{\prime }\right) }.
\end{equation}

In this case, we have
\begin{equation}\label{jose8}
\mathscr{L}\left[ \psi ^{\left( 1\right) }\left( \sigma ,t\right) \right] =s\left\{\psi \left( \sigma ,s\right) -
\left[ \varphi ^{\left( 1\right) }(t)\right] s^{-1}\psi \left( \sigma ,0\right)\right\} .
\end{equation}

As an example, will be check Eq.(22) of the paper \cite{yang}. Taking the Laplace transform on both sides of Eq.(\ref{jose}), we have
\begin{eqnarray}
\mathscr{L}\left[ ^{C}D_{a+}^{\xi \left( \sigma ,t\right) ;\varphi }\psi
\left( \sigma ,t\right) \right]  &=&L\left[ \frac{1}{\Gamma \left( 1-\xi
\left( \sigma ,t\right) \right) }\int_{a}^{t}\left( \varphi \left( \sigma
,t\right) -\varphi \left( \sigma ,s\right) \right) ^{-\xi \left( \sigma
,t\right) }\psi ^{\left( 1\right) }\left( \sigma ,s\right) ds\right]   \notag
\\
&=&\mathscr{L}\left[ \frac{\varphi \left( \sigma ,t\right) ^{-\xi \left(
\sigma ,t\right) }}{\Gamma \left( 1-\xi \left( \sigma ,t\right) \right) }%
\ast \psi ^{\left( 1\right) }\left( \sigma ,s\right) \right]   \notag \\
&=&\mathscr{L}\left[ \frac{\varphi \left( \sigma ,t\right) ^{1-\xi \left(
\sigma ,t\right) }}{\Gamma \left( -\xi \left( \sigma ,t\right) \right) }\ast
\psi ^{\left( 1\right) }\left( \sigma ,s\right) \right]   \notag \\
&=&\underset{\rm (I)}{\underbrace{\mathscr{L}\left[ \frac{\varphi \left( \sigma
,t\right) ^{1-\xi \left( \sigma ,t\right) }}{\Gamma \left( -\xi \left(
\sigma ,t\right) \right) }\right] }}\times \underset{\rm (II)}{\underbrace{%
\mathscr{L}\left[ \psi ^{\left( 1\right) }\left( \sigma ,s\right) \right] }},
\end{eqnarray}
where $\ast$ denotes the convolution.

Note that, we can use Eq.(\ref{jose8}), in part $\rm (II)$ of the product between $\rm (I)$ and $\rm (II)$. 
But, we can't use Eq.(15) as the authors report in the paper \cite{yang}, because as seen above, 
the definition of differential operator of variable order by means of the inverse Laplace transform, was used in a wrong way. 

In this sense, we conclude that
\begin{equation*}
\mathscr{L}\left[ ^{C}D_{a+}^{\xi \left( \sigma ,t\right) ;\varphi }\psi \left( \sigma ,t\right) \right] \neq \left[ \varphi ^{\left( 1\right) }\left( t\right) \right] s^{\xi \left( \sigma ,t\right) }\psi \left( s\right) -\left[ \varphi ^{\left( 1\right) }\left( t\right) \right] s^{\xi \left( \sigma ,t\right) -1}\psi \left( 0\right).
\end{equation*}

%%%%%%%%%%%%%%%%%%%%%%%%%%%%%%%%%%%%%%%%%%%%%%%%%%%%%%%%%%%%%%%%%%%%%%%%%%%%%%%%%%%%%%%%%%%%%%%%%%%%%%%%%%%%%%%%%%%%%%%%%%%%%%%%%%%%%%%%%%%%%%%%%%%%%%%%%%%%%%%%%%%%%%%%%%%%%%%%%%%%%%%%%%%%%%%%%%%%%%%%%%%%
\section{Concluding remarks}
Using the variable-order differential operator introduced by means of the inverse Laplace transform, we present a proof that the 
Laplace transform of $\varphi$-Caputo fractional derivative of the function $\psi\left(\sigma,t \right)$ of two-variable 
order $\xi \left (\sigma, t \right)$ $ \left (0 <\xi \left (\sigma, t \right) <1 \right) $ with respect to another function 
$\varphi$, is not the identity Eq.(22) of the mentioned paper \cite{yang}.

\bibliography{ref}
\bibliographystyle{plain}

\end{document}